\documentclass[journal,letterpaper,twoside, onecolumn]{IEEEtran}

\usepackage[cmex10]{amsmath}
\usepackage{fixltx2e}
\usepackage{url}
\usepackage[T1]{fontenc}
\usepackage{amssymb}
\usepackage{dsfont}
\usepackage{mathtools}
\usepackage{longtable}
\usepackage{booktabs}
\usepackage{tabularx}
\usepackage{ltxtable}
\usepackage{dcolumn}
\usepackage{lscape}
\usepackage{bm}

\newcommand{\R}{\mathds{R}}

\newcommand{\Z}{\mathds{Z}}
\newcommand{\Q}{\mathds{Q}}
\renewcommand{\O}{\mathcal{O}}
\renewcommand{\l}{$\ell_{1}$}

\newcommand{\st}{{\rm s.t.}}
\newcommand{\norm}[1]{\lVert {#1} \rVert}
\newcommand{\abs}[1]{\lvert {#1} \rvert}
\newcommand{\suchthat}{\; : \;}
\newcommand{\define}{\coloneqq}

\newcommand{\NP}{\textsf{NP}}
\newcommand{\coNP}{\textsf{coNP}}

\renewcommand{\P}{\textsf{P}}
\newcommand{\ceil}[1]{\lceil{#1}\rceil}

\newcommand{\A}{\bm{A}}

\newcommand{\B}{\bm{B}}
\renewcommand{\H}{\bm{H}}
\newcommand{\I}{\bm{I}}
\newcommand{\D}{\bm{D}}
\renewcommand{\L}{\bm{L}}
\newcommand{\E}{\bm{E}}
\newcommand{\x}{\bm{x}}
\renewcommand{\b}{\bm{b}}

\newtheorem{theorem}{Theorem}
\newtheorem{lemma}{Lemma}
\newtheorem{corollary}{Corollary}
\newtheorem{remark}{Remark}
\newtheorem{example}{Example}
\newtheorem{proposition}{Proposition}

\DeclareMathOperator{\spark}{spark}
\DeclareMathOperator{\rank}{rank}
\DeclareMathOperator{\diag}{Diag}
\DeclareMathOperator{\poly}{poly}

\begin{document}
\title{The Computational Complexity of 
the Restricted Isometry Property, the Nullspace Property, and Related Concepts in Compressed Sensing}
\author{Andreas~M.~Tillmann\IEEEmembership{} and
  Marc~E.~Pfetsch~\IEEEmembership{}
  \thanks{Manuscript received June 28, 2012; 
    last revised November 2, 2013; accepted November 3, 2013. This
    work has been partially funded by the Deutsche
    Forschungsgemeinschaft (DFG) within the project ``Sparse Exact and
    Approximate Recovery'' under grant PF 709/1-1. Parts of this work
    have been presented at the SPARS workshop, July 8-11, 2013, in
    Lausanne, Switzerland.}%
  \thanks{A. M. Tillmann is with the Research Group Optimization
    at TU Darmstadt, Dolivostr. 15, 64293 Darmstadt, Germany (phone: +49-6151-16-75833; e-mail:
    tillmann@mathematik.tu-darmstadt.de).}
  \thanks{M. E. Pfetsch is also with the Research Group Optimization at TU
    Darmstadt (e-mail:
    pfetsch@mathematik.tu-darmstadt.de).}
  \thanks{
    This work has been accepted for publication by the IEEE. Copyright
    will be transferred without notice, after which this version may
    no longer be accessible.}%
}

\markboth
{The Computational Complexity of RIP, NSP, and Related Concepts in Compressed Sensing}
{The Computational Complexity of RIP, NSP, and Related Concepts in Compressed Sensing}


\maketitle

\begin{abstract}
  \boldmath This paper deals with the computational complexity of
  conditions which guarantee that the \NP-hard problem of finding the
  sparsest solution to an underdetermined linear system can be solved
  by efficient algorithms. In the literature, several such conditions
  have been introduced. The most well-known ones are the mutual
  coherence, the restricted isometry property (RIP), and the nullspace
  property (NSP). While evaluating the mutual coherence of a given
  matrix is easy, it has been suspected for some time that evaluating
  RIP and NSP is computationally intractable in general. We confirm
  these conjectures by showing that for a given matrix~$\A$ and
  positive integer~$k$, computing the best constants for which the RIP
  or NSP hold is, in general, \NP-hard. These results are based on the
  fact that determining the spark of a matrix is \NP-hard, which is
  also established in this paper. Furthermore, we also give several
  complexity statements about problems related to the above concepts.
\end{abstract}


%
\IEEEpeerreviewmaketitle

\section{Introduction}\label{sect:intro}

\IEEEPARstart{A}{central} problem in compressed sensing (CS), see,
e.g., \cite{D06,E10,FR11}, is the task of finding a sparsest solution
to an underdetermined linear system, i.e.,
\begin{equation}\label{prob:P0}
  \min\,\norm{\x}_0\quad\st\quad \A\x=\b,\tag{${\rm P}_0$}
\end{equation}
for a given matrix $\A \in \R^{m \times n}$ with $m \leq n$,
where~$\norm{\x}_0$ denotes the $\ell_0$-quasi-norm, i.e., the number
of nonzero entries in~$\x$. This problem is well-known to be \NP-hard,
cf.~[MP5] in~\cite{GJ79}; the same is true for the \emph{denoising}
variant where the right hand side $\b$ is assumed to be contaminated
by noise and one employs the constraint
$\norm{\A\x-\b}_2\leq\varepsilon$ instead of $\A\x=\b$,
see~\cite{N95}.  Thus, in practice, one often resorts to efficient
heuristics. One of the most popular approaches is known as \emph{basis
  pursuit} (BP) or \l-minimization, see, e.g., \cite{CDS98}, where
instead of~\eqref{prob:P0} one considers
\begin{equation}\label{prob:P1}
  \min\,\norm{\x}_1\quad\st\quad \A\x=\b \tag{${\rm P}_1$}.
\end{equation}
Here, $\norm{\x}_1 = \sum_{i=1}^n \abs{x_i}$ denotes the
$\ell_1$-norm. It can be shown that under certain conditions, the
optimal solutions of~\eqref{prob:P0} and~\eqref{prob:P1} are unique
and coincide. In this case, one says that
\emph{$\ell_0$-$\ell_1$-equivalence} holds or that the
$\ell_0$-solution can be \emph{recovered} by the $\ell_1$-solution.
Similar results exist for the sparse approximation (denoising)
version with constraint $\norm{\A\x-\b}_2\leq\varepsilon$.

Many such conditions employ the famous \emph{restricted isometry
  property (RIP)} (see~\cite{CT05} and also~\cite{BDDVW08,BCT11}),
which is satisfied with order $k$ and a constant $\delta_k$ by a given
matrix $\A$ if
\begin{equation}\label{eq:rip}
  (1-\delta_k) \norm{\x}_2^2 ~\leq~ \norm{\A\x}_2^2 ~\leq~ (1+\delta_k)
  \norm{\x}_2^2 
\end{equation}
holds for all $\x$ with $1 \leq \norm{\x}_0 \leq k$. One is usually
interested in the smallest possible constant~$\underline{\delta}_k$,
the \emph{restricted isometry constant (RIC)}, such
that~\eqref{eq:rip} is fulfilled. For instance, a popular result
states that if $\underline{\delta}_{2k} < \sqrt{2}-1$, all $\x$ with
at most $k$ nonzero entries can be recovered (from $\A$ and $\b\define
\A\x$) via basis pursuit, see, e.g., \cite{C08}. A series of papers
has been devoted to developing conditions of this flavor, often by
probabilistic analyses. One of the most recent results,
see~\cite{CWX10}, shows that $\ell_0$-$\ell_1$-equivalence for
$k$-sparse solutions already holds if $\delta_k < 0.307$ is feasible
(note that the RIP order is $k$, not $2k$ as above). Several
probabilistic results show that certain random matrices are highly
likely to satisfy the RIP with desirable values of
$\underline{\delta}_k$, see, for instance, \cite{BT10,BCT11}. There
also has been work on deterministic matrix constructions aiming at
relatively good RIPs, see, e.g., \cite{DV07,I09,BDFKK11}. The RIP also
provides sparse recovery guarantees for other heuristics such as
Orthogonal Matching Pursuit and variants \cite{DW10,T04}, as well as
in the denoising case, see, e.g.,~\cite{C08}.

In the literature, it is often mentioned that evaluating the RIP,
i.e., computing the constant $\underline{\delta}_k$ for some $\A$ and
$k$, is presumably a computationally hard problem. Most papers seem to
refer to \NP-hardness, but this is often not explicitly stated. This
motivated the development of several (polynomial-time) approximation
algorithms for~$\underline{\delta}_k$, e.g., the semidefinite
relaxations in \cite{dAEGJL07,LB08}. However, while a widely accepted
conjecture in the CS community, \NP-hardness has, to the best of our
knowledge, not been proven so far.

Recently, some first results in this direction have been obtained:
In~\cite{KZ11} and \cite{KZ12}, hardness and non-approximability
results about the RIP were derived under certain (non-standard)
complexity assumptions; see also~\cite{BR13}. In work independent from
the present paper, \cite{BDMS12} shows that it is NP-hard to
verify~\eqref{eq:rip} for given~$\A$, $k$ and $\delta_k\in (0,1)$.

Another popular tool for guaranteeing $\ell_0$-$\ell_1$-equivalence is
the \emph{nullspace property (NSP)}, see, e.g.,
\cite{DH01,Zh08,CDDV09}, which characterizes recoverability
by~\eqref{prob:P1} (in fact, via $\ell_p$-minimization with $0 < p\leq
1$, see \cite{GN03,DG09}) for sufficiently sparse solutions
of~\eqref{prob:P0}. The NSP of order~$k$ is satisfied with
constant~$\alpha_k$ if for all vectors~$\x$ in the nullspace of~$\A$
(i.e., $\A\x=0$), it holds that
\begin{equation}\label{eq:NSPineq}
  \norm{\x}_{k,1} ~\leq~ \alpha_k \norm{\x}_1,
\end{equation}
where $\norm{\x}_{k,1}$ denotes the sum of the $k$ largest absolute
values of entries in $\x$. The NSP guarantees exact recovery of
$k$-sparse solutions to~\eqref{prob:P0} by solving~\eqref{prob:P1}
whenever~\eqref{eq:NSPineq} holds with some constant $\alpha_k <1/2$.

Similar to the RIP case, one is interested in the smallest
constant~$\underline{\alpha}_k$, the \emph{nullspace constant (NSC)},
such that~\eqref{eq:NSPineq} is fulfilled. Indeed, if and only if
$\underline{\alpha}_k < 1/2$, \eqref{prob:P1} with $\b\define\A
\tilde{\x}$, $\norm{\tilde{\x}}_0\leq k$, has the unique solution
$\tilde{\x}$, which coincides with that of~\eqref{prob:P0}; see, e.g.,
\cite{DH01,Zh08,JN11}. (Thus, the NSP provides a both necessary and
sufficient condition for sparse recovery, whereas the RIP is only
sufficient.) Moreover, error bounds for recovery in the denoising case
can be given, see, e.g., \cite{CDDV09}.

Again, the computation of~$\underline{\alpha}_k$ is suspected to be
\NP-hard, and several heuristics have been developed to compute good
bounds on $\underline{\alpha}_k$, e.g., the semidefinite programming
approaches in~\cite{dABEG08,dAEG11}, or an LP-based relaxation
in~\cite{JN11}.  However, as far as we know, no rigorous proof of
(\NP-)hardness has been given.

In this paper, we show that it is \NP-hard to compute the RIC and NSC
of a given matrix with given~$k$; see Sections~\ref{sect:ripNP}
and~\ref{sect:nspNP}, respectively. More precisely, we show that
unless \P$ = $\NP, there is no polynomial time algorithm that
computes~$\underline{\delta}_k$ or~$\underline{\alpha}_k$ for
\emph{all} given instances $(\A,k)$. We also prove that
\emph{certifying} the RIP given $\A$, $k$ \emph{and} some $\delta_k\in
(0,1)$ is \NP-hard.

Prior to this, in Section~\ref{sect:sparkNP}, we prove \NP-hardness of
computing the \emph{spark} of a matrix, i.e., the smallest number of
linearly dependent columns. In fact, our main results concerning the
complexity of determining the RIC or NSC follow from reduction of a
decision problem concerning the existence of small linearly dependent
column subsets. The term spark was first defined in~\cite{DE03}, where
strong results considering uniqueness of solutions to~\eqref{prob:P0}
were proven. Ever since, its value has been claimed to be \NP-hard to
calculate, but, to the best of our knowledge, without a proof or
reference for this fact. It seems to have escaped researchers' notice
that~\cite{K95} contains a proof that deciding whether the spark
equals the number of rows is 
\NP-hard, by reduction from the Subset Sum Problem (cf.~[MP9]
in~\cite{GJ79}). Moreover, \cite{CFGK03} provides a different proof
for this special case, by a reduction from the (homogeneous) Maximum
Feasible Subsystem problem \cite{AK95}. 
Even earlier, in~\cite{ColP86}, the authors claim to have a proof, but
give credit to the dissertation~\cite{McC83} for establishing
\NP-hardness of spark computations. However, after closer inspection,
the result in~\cite{McC83} is in fact \emph{not} about the spark, but
the girth of so-called transversal matroids of bipartite graphs. Only
recently, a variant of the latter proof has resurfaced
in~\cite{ACM11}, where it is used to derive (non-deterministic)
complexity results for constructing so-called \emph{full spark
  frames}, i.e., matrices exhibiting the highest possible spark. Every
transversal matroid can be represented by a matrix over any infinite
field or finite field with sufficiently large cardinality~\cite{PW70},
but there is no known deterministic way to construct such a matrix.

We adapt the proof idea from~\cite{McC83}, a reduction from the
$k$-Clique Problem, to vector matroids and thus establish that spark
computation is \NP-hard (without the restriction that the spark equals
the row size). Our proof also makes use of results from~\cite{CFGK03},
see Section~\ref{sect:sparkNP} for more details.

Moreover, we gather several more complexity statements regarding
problems related to the spark or RIP in Sections~\ref{sect:sparkNP}
and~\ref{sect:ripNP}, some of which are apparently new as well. In
particular, we also show that solving the sparse principal component
analysis problem (see, e.g., \cite{dAEGJL07,dAEG11,LT13}) is strongly
\NP-hard, which is another widely accepted statement that appears to
be lacking rigorous proof so far, and we extend this proof to show
that the \NP-hardness of RIC computation in fact holds in the strong
sense; see Section~\ref{subsect:asymmRIP}. Recall that \emph{strong}
\NP-hardness implies that (unless \P$=$\NP) there cannot exist a fully
polynomial-time approximation scheme (FPTAS), i.e., an algorithm that
solves a minimization problem within a factor of~$(1 + \varepsilon)$
of the optimial value in polynomial time with respect to the input
size and $1/\varepsilon$, see \cite{GJ79}; an FPTAS often exists for
weakly \NP-hard problems. Strong \NP-hardness can also be
understood as an indication that a problem's intractability does not
depend on ill-conditioning (due to the occurence of very large
numbers) of the input data.

Throughout the article, for an $m\times n$ matrix $\A$ and a subset $S
\subseteq \{1, \dots, n\}$, we denote by $\A_S$ the submatrix of~$\A$
formed by the columns indexed by~$S$. Sometimes we additionally
restrict the rows to some index set $R$ and write $\A_{RS}$ for the
resulting submatrix. Similarly, $\x_S$ denotes the part of a vector
$\x$ containing the entries indexed by~$S$. By $\A^\top$ and
$\x^\top$, we denote the transpose of a matrix $\A$ or vector $\x$,
respectively. For graph theoretic concepts and notation we refer
to~\cite{KorV08}, for complexity theory to~\cite{GJ79}, and for
matroid theory to~\cite{O92}.

\section{Complexity issues related to the spark}\label{sect:sparkNP}

In this section, we deal with complexity issues related to linearly
dependent columns of a given matrix~$\A\in\Q^{m \times n}$.
Inclusion-wise minimal collections of linearly dependent columns are
called circuits. More precisely, a \emph{circuit} is a set $C
\subseteq \{1, \dots, n\}$ of column indices such that $\A_C \x = 0$
has a nonzero solution, but every proper subset of~$C$ does not have
this property, i.e.,
$\rank(\A_C)=\abs{C}-1=\rank(\A_{C\setminus\{j\}})$ for every $j\in
C$. For notational simplicity, we will sometimes identify circuits~$C$
with the associated solutions~$\x \in \R^n$ of $\A \x = 0$ having
support~$C$. The \emph{spark} of~$\A$ is the size of its smallest
circuit.
\begin{example}\label{ex:spark}
  Consider the matrix
  \[
  \A = \left(\begin{array}{rrrr} 1 & 1 & 0 & 0\\ 1 & 1 & 0 & 1\\ 0 & 0 & 1 & 1\end{array}\right).
  \]
  Clearly, the first two columns yield the minimum-size circuit (in
  fact, the only one), i.e., $\spark(\A)=2$. In particular, note that
  generally, $\spark(\A)\leq k$ does \emph{not} guarantee that there
  also exists a \emph{vector with $k$ nonzeros} in the nullspace of
  $\A$; e.g., take $k=3$ for the above $\A$. On the other hand, it is
  immediately clear that a nullspace vector with support size $k$ does
  not yield $\spark(\A)=k$, but only $\spark(\A)\leq k$. This
  distinction between circuits and nullspace vectors in general will
  be crucial in the proofs below.
\end{example}

The main result of this section is the following.

\begin{theorem}\label{thm:circuitNP}
  Given a matrix $\A \in \Q^{m \times n}$ and a positive integer $k$,
  the problem to decide whether there exist a circuit of~$\A$ of size
  at most $k$ is \NP-complete.
\end{theorem}

For our proof, we employ several auxiliary results:
\begin{lemma}\label{lem:vN76}
  The vertex-edge incidence matrix of an undirected simple graph with
  $N$ vertices, $B$ bipartite components, and~$Q$ isolated vertices
  has rank~$N-B-Q$.
\end{lemma}

This result seems to be rediscovered every once in a while. The
earliest proof we are aware of is due to van Nuffelen~\cite{vN76} and
works through various case distinctions considering linear
dependencies of the rows and consequences of the existence of isolated
or bipartite components.

\begin{lemma}\label{lem:IncidenceRank}
  Let $G=(V,E)$ be a simple undirected
  graph 
  with vertex set~$V$ and edge set~$E$. Let $\A$ be its vertex-edge
  incidence matrix, and let $k>4$ be some integer. Suppose $G$ only
  has connected components with at least four vertices each,
  $\abs{E}=\binom{k}{2}$, and $\rank(\A)=k$. Then the graph~$G$ has
  exactly $\abs{V}=k$ vertices.
\end{lemma}

\begin{IEEEproof}
  Let $G = (V, E)$ and $k > 4$ be the graph and integer given in the
  statement of the lemma. Assume that $G$ has no component with less
  than four vertices, has $\abs{E} = \binom{k}{2}$ edges, and that its
  incidence matrix $\A$ 
  has $\rank(\A)=k$.

  Since consequently, $G$ has no isolated vertices, Lemma~\ref{lem:vN76} tells us
  that the number of vertices is
  \[
  N ~=~\rank(A)+B ~=~ k + B,
  \]
  where $B$ is the number of bipartite components in $G$. Assume that
  $B > 0$, since otherwise the lemma is trivially true.

  We claim that the number of edges in~$G$ can be at most
  \begin{equation}\label{eq:EdgeNum}
    \abs{E} ~\leq~ \binom{N}{2} - \frac{4(N-4)}{2}B - 2\,B.
  \end{equation}
  To see this, recall that $G$ can have at most $\binom{N}{2}$
  edges. Each connected component has at least four vertices.  Since
  there are no edges between such a component and vertices outside,
  the total number of possible edges is reduced by at least $4(N-4)/2$
  per component (the factor $1/2$ ensures that we do not count any
  edges twice). Since $G$ has at least~$B$ connected components, the
  possible number of edges is hence decreased at least by the second
  term in~\eqref{eq:EdgeNum}. Moreover, since each bipartite component
  has at least four vertices, at least two of the potential edges
  cannot be present inside each such component, which yields the last
  term in~\eqref{eq:EdgeNum}. Note that the bound~\eqref{eq:EdgeNum}
  is sharp if $G$ consists only of bipartite components with four
  vertices each.

  Expanding~\eqref{eq:EdgeNum} using $N=k+B$, we obtain
  \begin{align*}
    \abs{E} ~\leq~ \binom{N}{2} - \frac{4(N-4)}{2}B - 2\,B ~&=~ \frac{(k+B)(k+B-1)}{2} - \frac{4(k+B-4)}{2}B - 2\,B\\
    &=~ \tfrac{1}{2} k^2 - \tfrac{1}{2} k - k\, B - \tfrac{3}{2} B^2 +  \tfrac{11}{2}\, B ~=~ \binom{k}{2} - (k\, B + \tfrac{3}{2} B^2 -  \tfrac{11}{2}\, B),
  \end{align*}
  and observe that
  \[
  k\, B + \tfrac{3}{2} B^2 - \tfrac{11}{2}\, B ~\geq~ 5\, B +
  \tfrac{3}{2} B^2 - \tfrac{11}{2}\, B ~=~ \tfrac{3}{2} B^2 -
  \tfrac{1}{2}\, B ~>~ 0
  \]
  if $B > 0$. Thus, there are strictly less than $\binom{k}{2}$ edges,
  contradicting the requirement $\abs{E}=\binom{k}{2}$. Hence, $B=0$.
\end{IEEEproof}

\begin{lemma}\label{lem:CFGK03}
  Let $\H = (h_{ij}) \in\Z^{m\times n}$ be a full-rank integer matrix
  with $m\leq n$ and let $\alpha\define \max\, \abs{h_{ij}}$. Let
  $q\in\{m,\dots,n\}$ and define
  \[
  \H(x) \define \left(\begin{array}{ccccc}
      \multicolumn{5}{c}{\H} \\
    1 & x & (x)^2 & \dots & (x)^{n-1}\\
    1 & x+1 & (x+1)^2 & \dots & (x+1)^{n-1}\\
    \vdots & \vdots & \vdots & \dots & \vdots\\
    1 & x+q-m-1 & (x+q-m-1)^2 & \dots & (x+q-m-1)^{n-1}\end{array}\right) \in\Z^{q\times n}.
  \]
  For any column subset $S$ with $\abs{S}=q$, if $\rank(\H_S)=m$ and
  $\abs{x}\geq \alpha^m q^{qn}+1$, then $\H(x)_S\in\Z^{q\times q}$ has
  full rank $q$.
\end{lemma}

\begin{IEEEproof}
  This result is a combination of Lemma~1 and (ideas from the proof
  of) Proposition~4 from~\cite{CFGK03}.  For clarity, we give the
  details here. Consider some $S\subseteq\{1,\dots,n\}$ with
  $\abs{S}=q$ (w.l.o.g., $q>m\geq 1$; otherwise there is nothing to
  show). Assume that $\rank(\H_S)=m$ and note that the last $q-m$ rows
  of $\H(x)_S$ form a submatrix of a generalized Vandermonde matrix
  with distinct nodes (see, e.g., \cite{DM99}), which is easily seen
  to have full rank as well; see also~\cite{CFGK03}. Consider the
  polynomial $p(x)\define{\rm
    det}(\H(x)_S)$. From~\cite[Lemma~1]{CFGK03}, we know that there
  exists some $x$ for which the subspace spanned by the last $q-m$
  rows of $\H(x)_S$ and the row space of $\H_S$ are transversal, i.e.,
  they only intersect trivially. In particular, this shows that $p$
  cannot be identical to the zero polynomial (both transversal parts
  have full rank). Let $d$ be the degree of $p(x)$ (which depends on
  the choice of $S$); thus, $p(x)=\beta_0 +\beta_1 x+\dots+\beta_d
  x^d$ with $\beta_d\neq 0$. 
  Expanding the determinant ${\rm det}(\H(x)_S)$ using Leibniz's
  formula, one can derive that $\abs{\beta_i}\leq \alpha^m q^{qn}$ for
  all~$i$, by noting that the expansion consists of $q!<q^q$ summands
  which each are the product of precisely one entry per matrix row and
  column---in absolute value terms, we can extract a factor of
  $\alpha^m$ from this sum (from the $m$ rows corresponding to
  $\H_S$), and upper-bound all absolute values of coefficients of $x$
  by the highest possible value (occurring when the last $q-m$ columns
  of $\H(x)$ are contained in $S$), which can be no larger than
  $(q-m)^{(n-1)(q-m)}<q^{qn-q}$. Moreover, it is easy to see that
  $\beta_i\in\Z$ for all~$i$. Applying Cauchy's bound~\cite{HM97} to
  the monic polynomial obtained from dividing $p(x)$ by $\beta_d$
  yields
  \[
  \abs{r} ~<~ 1+\max_{0\leq i\leq d-1}\left\lvert\frac{\beta_i}{\beta_d}\right\rvert ~\leq~ 1+\max_{0\leq i\leq d-1}\abs{\beta_i} ~\leq~ 1+\alpha^m q^{qn}\qquad\text{for all }r\text{ with }p(r)=0.
  \]
  Thus, any $x$ with $\abs{x}\geq \alpha^m q^{qn}+1$ is not a root of
  $p(x)$; equivalently, $\rank(\H(x)_S)=q$ for such $x$.
\end{IEEEproof}

\begin{IEEEproof}[Proof of Theorem~\ref{thm:circuitNP}]
  The problem is clearly in \NP: Given a subset $C$ of column indices
  of~$\A$, it can be verified in polynomial time that $\abs{C}\leq k$
  and that $\rank(\A_C) = \abs{C}-1 = \rank(\A_{C\setminus\{j\}})$ for
  every $j \in C$ (by Gaussian elimination, see, e.g.,
  \cite{GroLS93}).

  To show hardness, we reduce the \NP-complete $k$-Clique Problem
  (cf.\ [GT19] in~\cite{GJ79}, or~\cite{K72}): Given a simple
  undirected graph~$G$, decide whether~$G$ has a \emph{clique}, i.e.,
  a vertex-induced complete subgraph, of size~$k$. We may assume
  without loss of generality that $k>4$.

  For the given graph~$G$ with~$n$ vertices and~$m$ edges construct a
  matrix~$\A=(a_{ie})$ of size $(n + \binom{k}{2} - k - 1) \times m$
  as follows: Index the first~$n$ rows of~$\A$ by the vertices of~$G$
  and its columns by the edges of~$G$ (we will also identify the
  vertices and edges with their indices). Let the first~$n$ rows of
  $\A$ contain the vertex-edge incidence matrix of~$G$ (i.e., set
  $a_{ie}=1$ if $i\in e$, and $0$ otherwise). For the non-vertex rows
  $n+i$, $i \in \{1,\dots,\binom{k}{2}-k-1\}$, set $a_{(n+i)e}
  =(U+i-1)^{e-1}$ with $U\define k^{2k^2 m}+1$; note that this
  corresponds to the bottom part of $A$ consisting of a Vandermonde
  matrix (each row consists of increasing powers of the distinct
  numbers $U,\dots,U+\binom{k}{2}-k-2$). Clearly, this matrix~$\A$ can
  be constructed in polynomial time, and its encoding length is
  polynomially related to that of the input (in particular, that of
  its largest entry is $\O(k^2m^2\log_2(k))$).

  We first show that $G$ has a $k$-clique if and only if $\A$ has a
  circuit of size~$\binom{k}{2}$. Suppose that~$G$ has a $k$-clique,
  $k > 4$, say on the vertices in the set~$R$ (so that $\abs{R} = k$),
  and with its $\binom{k}{2}$ edges in the set~$C$. Since $\A_C$ has
  all-zero rows for each vertex outside of~$R$, $\abs{R}+\text{(number
    of non-vertex rows)} = \binom{k}{2} - 1 =
  \abs{C}-1\geq\rank(\A_C)$. Clearly, a clique is never bipartite (it
  always contains odd cycles, for $k\geq 3$). Hence, by
  Lemma~\ref{lem:vN76}, the rows of $\A_C$ indexed by $R$ are linearly
  independent. Now observe that removing any edge from a $k$-clique
  does not affect the rank of the associated incidence matrix, since
  the subgraph remains connected and non-bipartite with less vertices
  than edges (for $k \geq 4$). Thus, by Lemma~\ref{lem:vN76}, the rank
  of the nonzero vertex row part of $\A_C$ remains $k$ if any column
  from $C$ is removed. Therefore, since $a_{ie}\leq 1$ for all $i\leq
  n$ and all $e\leq m$, Lemma~\ref{lem:CFGK03} applies to
  $\A_{C\setminus\{e\}}$ for every $e\in C$ (with $x=U$, $\H(x)=\A$,
  $S=C\setminus\{e\}$, and $q = \binom{k}{2} - 1$) and yields
  $\rank(\A_{C\setminus\{e\}})=q=\abs{C}-1$, 
  whence also $\rank(\A_C)=\abs{C}-1$. Thus, $C$ is a circuit.

  Conversely, suppose that~$\A$ has a circuit~$C$ of size
  $\abs{C}=\binom{k}{2}$ with $k > 4$. Then, by definition of a
  circuit, $\rank(\A_C)=\abs{C}-1$, so $\A_C$ has at least $\abs{C} -
  1$ nonzero rows. Since these include the $\abs{C} - k-1$ non-vertex
  rows, the set~$R$ of nonzero vertex rows of~$\A_C$ has size
  $\abs{R}\geq \big(\abs{C}-1\big) - \big(\abs{C} - k-1\big) = k$.
  Let $\A_{RC}$ and $\A_{NC}$ denote the vertex and non-vertex row
  submatrices of $\A_C$, respectively. Since
  $\rank(\A_C)\leq\rank(\A_{NC})+\rank(\A_{RC})=\abs{C}-k-1+\rank(\A_{RC})$
  and $\abs{R}\geq k$, clearly $\rank(\A_{RC})\geq k$. Suppose that
  $\rank(\A_{RC})>k$; then there would exist a subset $R'\subseteq R$
  with $\abs{R'}=k+1$ and $\rank(\A_{R'C})=k+1$. But by
  Lemma~\ref{lem:CFGK03}, the square matrix
  $(\A_{R'C}^\top,\A_{NC}^\top)^\top$ would then have full rank
  $\abs{C}$, whence $\rank(\A_C)=\abs{C}>\abs{C}-1$, contradicting the
  fact that $C$ is a circuit. Thus, the upper part $\A_{RC}$ of~$\A_C$
  must in fact have rank \emph{exactly}~$k$.

  Observe that the subgraph $(R,C)$ of~$G$ with vertex set~$R$ and
  edge set~$C$ cannot contain components with less than $4$ vertices:
  such a subgraph $(R',C')$ could have at most as many edges as
  vertices, so that the associated incidence matrix $\A'_{C'}$ has
  full \emph{column} rank. Removing a column corresponding to an edge
  $e \in C'$ would reduce the rank, i.e., $\rank(\A'_{C' \setminus
    \{e\}}) < \rank(\A'_{C'})$. Moreover, note that $\A_{RC}$ has
  block diagonal form where the blocks are the incidence matrices of
  the separate graph components within $(R,C)$, so that the rank is
  the sum of the ranks of the blocks (one of which is $\A_{R'C'}$). In
  particular, the non-vertex row part of $\A_C$ maintains full (row)
  rank when any column is removed from $C$, so that deleting $e\in C'$
  would yield $\rank(\A_{C\setminus\{e\}}) = \rank(\A_C)-1$,
  contradicting the fact that $C$ is a circuit. Thus,
  Lemma~\ref{lem:IncidenceRank} applies to the graph $(R,C)$ and
  yields that $\abs{R} = k$. This implies that the vertices in~$R$
  form a $k$-clique, because~$R$ can induce at most~$\binom{k}{2}$
  edges and the~$\binom{k}{2}$ edges in~$C$ are among them.

  We now show that each circuit of $\A$ has size \emph{at least}
  $\binom{k}{2}$. This proves the claim, since by the arguments above,
  it shows that there exists a circuit of size \emph{at most} (in
  fact, exactly) $\binom{k}{2}$ if and only if $G$ has a $k$-clique,
  i.e., for the given construction a solution to the spark problem
  yields a solution to the clique problem as well.

  Suppose that~$\A$ has a circuit~$C$ with $c \define \abs{C} <
  \binom{k}{2}$. Let $d\define\binom{k}{2}-c > 0$. Clearly, not all
  vertex rows restricted to any column subset can be zero, and any
  submatrix of the non-vertex part of $\A$ with fewer than
  $\binom{k}{2}-k$ columns is of full (column) rank. Therefore,
  $c>\binom{k}{2}-k$ necessarily. Since~$C$ is a circuit, $\A_C$ has
  $c-1$ nonzero rows (similar to the arguments above, it can be seen
  that the lower bound $c-1$ holds with equality). Because the
  $\binom{k}{2}-k-1$ non-vertex rows are among these, and by
  Lemmas~\ref{lem:IncidenceRank} and~\ref{lem:CFGK03}, $\A_C$ has
  $(c-1)-\big(\binom{k}{2}-k-1\big) = k-d>0$ nonzero vertex
  rows. Denote the set of such rows by~$R$, and let~$r$ be the number
  of edges in the subgraph of~$G$ induced by the vertices
  in~$R$. Since $\abs{R} = k - d$ vertices can induce at most
  $\binom{k-d}{2}$ edges, $r \leq \binom{k-d}{2}$. But surely all the
  edges in~$C$ are among those induced by~$R$, so that $r\geq c =
  \binom{k}{2} - d$. Putting these inequalities together yields
  $\binom{k}{2} - d \leq r \leq \binom{k-d}{2}$. However, since we
  assumed $k > 4$, it holds that $\binom{k}{2}-d > \binom{k-d}{2}$,
  yielding a contradiction. Consequently, every circuit $C$ of $A$
  must satisfy $\abs{C} \geq \binom{k}{2}$.
\end{IEEEproof}

The smallest size (cardinality) of a circuit can be expressed as
\begin{equation}\label{eq:sparkdef}
  \spark(\A) ~\define~ \min\{\,\norm{\x}_0 \suchthat \A\x = 0,\ \x \neq 0\,\}.
\end{equation}
%
Clearly, there exists a circuit of size at most~$k$ if and only if the spark
is at most~$k$. Hence, Theorem~\ref{thm:circuitNP} immediately yields the
following.

\begin{corollary}\label{cor:sparkNP}
  Computing $\spark(\A)$ is \NP-hard.
\end{corollary}

\begin{remark}\label{rem:spark}
  The idea of reducing from the clique problem to prove
  Theorem~\ref{thm:circuitNP} is due to Larry Stockmeyer and appears
  in Theorem~3.3.6 of~\cite{McC83} (see
  also~\cite{ColP86,ACM11}). However, \cite{McC83} uses \emph{generic}
  matrices that, in fact, represent transversal matroids (of bipartite
  graphs) and therefore have certain properties needed in the
  proof. The entries of these generic matrices are not specified, and
  to date there is no known deterministic way to do so such that the
  matrix represents a transversal matroid. We replaced the
  corresponding machinery by our explicit matrix construction and the
  arguments using Lemmas~\ref{lem:vN76}, \ref{lem:IncidenceRank}
  and~\ref{lem:CFGK03} to become independent of transversal matroid
  representations and work directly on vector matroids.  Note that the
  proof of Theorem~\ref{thm:circuitNP} also shows \NP-completeness of
  the problem to decide whether $\A$ has a circuit $C$ with
  \emph{equal to} $k$; see also~Example~\ref{ex:spark}.
\end{remark}

\begin{remark}\label{rem:spark2}
  The results above are related to, but different from, the following.
  \begin{enumerate} 
  \item Theorem~1 in~\cite{K95} shows that for an $m\times n$ matrix
    $\A$, it is \NP-complete to decide whether $\A$ has an $m\times m$
    submatrix with zero determinant. This implies \NP-completeness of
    deciding whether $\spark(\A)\leq k$ for the special case
    $k=m$. This restriction of $k$ to the row number $m$ of $\A$ could
    in principle be removed by appending all-zero rows, but one would
    then no longer be in the interesting case where the matrix has
    full (row) rank.  Our proof admits spark values other than the row
    number for full-rank matrices; however, the row and column numbers
    in the reduction depend on the instance.

  \item In contrast to the results above, for graphic matroids, the
    girth can be computed in polynomial time~\cite{IR78,McC83}.

  \item The paper~\cite{V97} proves \NP-hardness of computing the
    girth of the binary matroid, i.e., a vector matroid over
    $\mathds{F}_2$. This, however, does not imply \NP-hardness over
    the field of rational or real numbers, and the proof cannot be
    extended accordingly. Similarly, in~\cite{BMvT78} it was shown
    that, over $\mathds{F}_2$, it is \NP-complete to decide whether
    there exists a vector with~$k$ nonzero entries in the nullspace of
    a matrix. However, while the proof for this result can
    straightforwardly be extended to the rational case, it does not
    imply hardness of computing the spark either: Since
    in~\cite{BMvT78}, there is no lower bound on the spark (such as we
    provide in the last paragraph of the proof of
    Theorem~\ref{thm:circuitNP}), a situation as in
    Example~\ref{ex:spark} is not explicitly avoided there.  (Note
    also that it was already remarked in~\cite{BMvT78} that the
    problem to decide whether an ($\mathds{F}_2$-)nullspace vector
    with \emph{at most}~$k$ nonzeros exists is not covered by their
    proof.)
  \end{enumerate}
\end{remark}

We also have the following result, which shows another relation
between minimum cardinality circuits and the task of finding sparsest
solutions to underdetermined linear systems.

\begin{corollary}\label{cor:sparkP0}
  Given a matrix $\bm{B}$, a specific column $\b$ of $\bm{B}$, and a
  positive integer~$k$, the problem of deciding whether there exists a
  circuit of $\bm{B}$ of size $k$ which contains $\b$ is \NP-complete
  in the strong sense. Consequently, it is strongly \NP-hard to
  determine the minimum cardinality of circuits that contain a
  specific column $\b$ of $\bm{B}$.
\end{corollary}

\begin{IEEEproof}
  Denote by $\A$ the matrix $\bm{B}$ without the column $\b$. Then it
  is easy to see that $\bm{B}$ has a circuit of size $k$ that
  contains~$\b$ if and only if $\A\x=\b$ has a solution with $k-1$
  nonzero entries. Deciding the latter is well-known to be
  \NP-complete in the strong sense (it amounts to the decision version
  of~\eqref{prob:P0}), see [MP5] in~\cite{GJ79}.
\end{IEEEproof}

\begin{remark}
  As mentioned in~\cite{DE03}, one can reduce spark computations
  to~\eqref{prob:P0} as follows: For each column of~$\A \in \Q^{m
    \times n}$ in turn, add a new row with a~$1$ in this column
  and~$0$ elsewhere. The right hand side~$\b$ is the $(m+1)$-th unit
  vector. Now solve each such~\eqref{prob:P0} problem, and take the
  solution with smallest support. (Note that this is a
  (Turing-)reduction, using Theorem~\ref{thm:circuitNP} to show
  \NP-hardness of~\eqref{prob:P0}.) Interestingly, we do not know an
  easy reduction for the reverse direction.
\end{remark}

Let us now briefly consider full spark frames. An $m\times n$ matrix
$\A$ with full rank $m$ ($m\leq n$) is said to be \emph{full spark} if
$\spark(\A)=m+1$, i.e., every submatrix consisting of at most $m$
columns of $\A$ has full rank. In~\cite{ACM11}, it was shown that
testing a matrix for this property is hard for {\NP} under randomized
reductions.  In fact, the following stronger result holds:
\begin{corollary}\label{cor:fullspark}
  Given a rational matrix $\A$, deciding whether $\A$ is a full spark
  frame is \coNP-complete.
\end{corollary}
\begin{IEEEproof}
  We can assume w.l.o.g. that $\A\in\Q^{m\times n}$ with rank $m\leq
  n$. Thus, $\A$ is full spark if $\spark(\A)=m+1$. Clearly,
  $\spark(\A)=m+1$ holds if and only if the question whether $\A$ has
  a singular $m\times m$ submatrix has a negative answer. Since the
  latter decision problem is \NP-complete by~\cite{K95}
  (or~\cite[Proposition~4]{CFGK03} and the results in~\cite{AK95}),
  deciding whether $\A$ is a full spark frame is \NP-hard.  Moreover,
  this problem is contained in \coNP, since the ``no'' answer can be
  certified in polynomial time by specifying a singular (square)
  submatrix and because singularity can be verified in polynomial
  time.
\end{IEEEproof}

Note that for the above proof, we cannot employ
Theorem~\ref{thm:circuitNP}, because this would require considering
the matrix construction used in the reduction from the clique problem
(see the proof of the theorem) with $k=n-1$, and the clique problem
can be solved in polynomial time for any $k=n-\ell$ with \emph{constant}
$\ell$.

\begin{remark}\label{rem:coNP}
  Above, and in the related complexity results to follow, we show that
  the decision problem under consideration has a \emph{negative}
  answer if and only if a known \NP-complete problem has a
  \emph{positive} answer. Since, by definition, the complementary
  problem of an \NP-complete problem is \coNP-complete, the respective
  hardness results follow---see also~\cite{GJ79}. Both \NP- and
  \coNP-completeness imply that no polynomial time algorithm exists
  unless \P$=$\NP\ (or equivalently \P$=$\coNP). Since a problem is
  \NP-hard if and only if it is \coNP-hard (every problem in \coNP\
  can be Turing-reduced to it; cf~\cite[Ch.~15.7]{KorV08}), we use the
  term \NP-hard throughout the article.
\end{remark}

In the next sections we shall see how we can deduce \NP-hardness of
computing restricted isometry or nullspace constants from the above
results. (In particular, the extension to $k<m$ provided by
Theorem~\ref{thm:circuitNP}, cf. Remark~\ref{rem:spark2}, will allow
for making statements about RIP or NSP orders other than the row
number.)

\section{\NP-hardness of computing the restricted isometry constant}\label{sect:ripNP}

Recall that for a positive integer~$k$, $\A$ satisfies the RIP of
order~$k$ with constant~$\delta_k$ if~\eqref{eq:rip} holds. Given $\A$
and $k$, the smallest such constant is the RIC
$\underline{\delta}_k$. Note that~\eqref{eq:rip} holds for $\delta_k =
0$ if and only if~$\A$ is orthogonal, and that $\delta_k < 0$ is
impossible.

The suspicions about computational intractability of the RIP are based
on the observation that a brute-force method would have to inspect all
submatrices induced by column subsets of sizes up to $k$. Of course,
by itself, this does not generally rule out the possible existence of
an efficient algorithm. However, we show below that (given $\A$ and
$k$) deciding whether there exists a constant~$\delta_k < 1$ such
that~\eqref{eq:rip} holds is \coNP-complete; consequently, computing
the RIC is \NP-hard. Moreover, we show that RIP certification (for a
given $\delta_k\in (0,1)$) is \NP-hard. Thus, unless \P$=$\NP, there
exists no polynomial time algorithm for any of these problems
(cf. Remark~\ref{rem:coNP}).

We will need the following technical result.
\begin{lemma}\label{lem:scaleA}
  Let $\A = (a_{ij}) \in \Q^{m \times n}$ be a matrix. Define
  $
  \alpha \define \max \{ \abs{a_{ij}} \suchthat i \in \{1, \dots, m\},\; j \in
  \{1, \dots, n\}\}
  $ 
  and set $C \define 2^{\ceil{\log_2(\alpha\, \sqrt{mn})}}$. Then $\tilde{\A}
  \define \tfrac{1}{C} \A$ has encoding length polynomial in that of~$\A$, and satisfies
  \[
  \norm{\tilde{\A} \x}_2^2 ~\leq~ (1+\delta) \norm{\x}_2^2\qquad\text{for all~} \x\in\R^n\text{ and }\delta\geq 0.
  \]
\end{lemma}
\begin{IEEEproof}
  First, observe that the largest singular value of~$\A$,
  $\sigma_{\max}(\A)$, satisfies
  \[
    \sigma_{\max}(\A) ~=~ \norm{\A}_2 ~\leq~ \alpha \sqrt{mn}~\leq~ 2^{\ceil{\log_2(\alpha\, \sqrt{mn})}} ~=~ C.
  \]
  It follows that
  \[
  \norm{\tilde{\A} \x}_2^2 ~\leq~ \norm{\tfrac{1}{C}\A}_2^2\, \norm{\x}_2^2 ~\leq~ \tfrac{1}{\sigma_{\max}(\A)^2}\norm{\A}_2^2\norm{\x}_2^2 ~=~ \norm{\x}_2^2 ~\leq~ (1+\delta)\norm{\x}_2^2
  \]
  for any $\delta\geq 0$. Moreover, the encoding length of~$C$ and
  of~$\tilde{\A}$ is clearly polynomially bounded by that of~$\A$, $m$
  and $n$.
\end{IEEEproof}

By the singular value interlacing theorem (see, e.g., \cite{Q87}),
$\sigma_{\max}(\A)$ is an upper bound for the largest singular value
of every submatrix of $\A$. Thus, the above lemma essentially shows
that by scaling the matrix~$\A$, one can focus on the lower part of
the RIP~\eqref{eq:rip} (a similar argument has been derived
independently in~\cite{BDMS12}). This leads to the following
complexity result.

\begin{theorem}\label{thm:ripNP}
  Given a matrix~$\A \in \Q^{m \times n}$ and a positive integer~$k$,
  the problem to decide whether there exists some rational constant
  $\delta_k<1$ such that $\A$ satisfies the RIP of order~$k$ with
  constant $\delta_k$ is \coNP-complete.
\end{theorem}

\begin{IEEEproof}
  We first show that the problem is in \coNP. To certify the ``no''
  answer, it~suffices to consider a vector~$\tilde{\x}$ with $1 \leq
  \norm{\tilde{\x}}_0 \leq k$ that tightly satisfies~\eqref{eq:rip}
  for~$\delta_k = 1$. This implies $\A \tilde{\x} = 0$. Clearly, since
  $\tilde{\x}$ is contained in the nullspace of~$\A$, we can assume
  that $\tilde{\x}$ is rational with encoding length polynomially
  bounded by that of~$\A$. Then, we can verify
  $1\leq\norm{\tilde{\x}}_0\leq k$ and $\A\tilde{\x} = 0$ in
  polynomial time, which shows that the ``no'' answer can be certified
  in polynomial time.

  To show hardness, we reduce the problem to decide whether there
  exists a circuit of size at most~$k$, which is \NP-complete by
  Theorem~\ref{thm:circuitNP}. Consider the matrix~$\tilde{\A}$ as
  defined in Lemma~\ref{lem:scaleA}; note that the circuits of~$\A$
  and~$\tilde{\A}$ coincide, since nonzero scaling does not affect
  linear dependencies among columns. We claim that there exists a
  circuit~$\tilde{\x}$ with $1 \leq \norm{\tilde{\x}}_0 \leq k$ if and
  only if~$\tilde{\A}$ \emph{violates}~\eqref{eq:rip} for \emph{all}
  $\delta_k < 1$. Since deciding the former question is \NP-complete,
  this completes the proof.

  Clearly, if such an~$\tilde{\x} \neq 0$ exists, then
  \[
  (1 - \delta_k) \norm{\tilde{\x}}_2^2 ~\leq~ \norm{\tilde{\A} \tilde{\x}}_2^2 ~=~ 0
  \]
  implies that we must have $\delta_k \geq 1$.

  For the converse, assume that there does not exist~$\delta_k < 1$
  for which~\eqref{eq:rip} holds. By Lemma~\ref{lem:scaleA}, the upper
  part of~\eqref{eq:rip} is always satisfied. Consequently, there must
  exist a vector~$\hat{\x}$ with $1 \leq \norm{\hat{\x}}_0 \leq k$
  such that the lower part is tight, i.e.,
  \[
  0 ~\geq~ (1 - \delta_k) \norm{\hat{\x}}_2^2 ~=~ \norm{\tilde{\A} \hat{\x}}_2^2 ~\geq~ 0.
  \]
  This implies that $\tilde{\A} \hat{\x} = 0$. Thus, there also
  exists a \emph{circuit}~$\tilde{\x}$ with $1\leq \norm{\tilde{\x}}_0
  \leq k$, which shows the claim.
\end{IEEEproof}

Usually, one is interested in the smallest constant~$\delta_k$ for
which~\eqref{eq:rip} holds, i.e., the \emph{restricted isometry
  constant (RIC)}
\begin{equation}\label{eq:deltak}
  \underline{\delta}_k \define \min_{\delta \geq 0} ~\delta\quad\st\quad (1 - \delta) \norm{\x}_2^2 \leq \norm{\A\x}_2^2 \leq (1 + \delta)\norm{\x}_2^2\quad\text{for all } \x\text{ with } \norm{\x}_0 \leq k.
\end{equation}
Recall that $\delta < 0$ is impossible, resulting in the condition
$\delta \geq 0$ in~\eqref{eq:deltak}. We immediately obtain the
following complexity result.

\begin{corollary}\label{cor:ripNP01NPhard}
  For a given matrix $\A \in \Q^{m \times n}$ and positive integer $k$,
  it is \NP-hard to compute the RIC~$\underline{\delta}_k$.
\end{corollary}

\subsection{RIP certification}\label{subsect:RIPcertification}

In this section, we show that the RIP certification problem, i.e.,
deciding whether a matrix~$\A$ satisfies the RIP with given order~$k$
and \emph{given} constant $\delta_k \in (0,1)$, is (\textsf{co})\NP-hard.
The main arguments used in the proofs of the following Lemma and
Theorem have been independently derived by us and the authors
of~\cite{BDMS12}.

The following observation is essential.
\begin{lemma}\label{lem:RIPlowerBound}
  Given a matrix $\A\in\Q^{m\times n}$ and a positive integer $k \leq n$,
  if $\spark(\A)>k$, there exists a constant $\varepsilon>0$ with encoding
  length polynomially bounded by that of $\A$ such that
  $\norm{\A\x}_2^2 \geq \varepsilon\, \norm{\x}_2^2$ for all $x$ with
  $\norm{\x}_0\leq k$.
\end{lemma}
\begin{IEEEproof}
  Assume without loss of generality that $\A$ has only integer entries
  (this can always be achieved by scaling with the least common
  denominator 
  of the matrix entries, which influences~$\varepsilon$ by a
  polynomial factor only).
  
  Define $\alpha$ as in Lemma~\ref{lem:scaleA}. Note that
  $\spark(\A)>k$ implies that every submatrix $\A_S$ with $S \subseteq
  \{1, \dots, n\}$, $\abs{S}\leq k$, has linearly independent
  columns. Consider an arbitrary such~$S$. Then, $\A_S^\top \A_S$ is
  positive definite, so its smallest eigenvalue fulfills
  $\lambda_{\min}(\A_S^\top \A_S^{})>0$, and also $\det(\A_S^\top
  \A_S) > 0$. Moreover, since the absolute values of entries of~$\A$
  are integers in $\{0,1,\dots,\alpha\}$, the entries of $\A_S^\top
  \A_S^{}$ are also integral and lie in
  $\{0,1,\dots,m\,\alpha^2\}$. Therefore, it must in fact hold that
  $\det(\A_S^\top \A_S^{})\geq 1$ and $\lambda_{\max}(\A_S^\top
  \A_S^{})\geq 1$. It follows that
  \begin{align*}
    1 ~&\leq~ \det(\A_S^\top \A_S^{}) ~=~ \prod_{i=1}^{\abs{S}}\lambda_i(\A_S^\top \A_S^{})
    ~\leq~ \lambda_{\min}(\A_S^\top \A_S^{})\cdot\lambda_{\max}(\A_S^\top \A_S^{})^{k-1}\\
    &\leq~ \lambda_{\min}(\A_S^\top \A_S^{})\left(\abs{S}\cdot\max_{1\leq i,j\leq \abs{S}}\abs{(\A_S^\top \A_S^{})_{ij}}\right)^{k-1}
    ~\leq~ \lambda_{\min}(\A_S^\top \A_S^{})\left(k\, m\, \alpha^2\right)^{k-1}.
  \end{align*} 
  Consequently, we have that
  \begin{equation}\label{eq:RIPeps}
    \lambda_{\min}(\A_S^\top \A_S^{}) ~\geq~ \frac{1}{(k\, m\, \alpha^2)^{k-1}} ~\eqqcolon~ \varepsilon ~>~ 0.
  \end{equation}
  
  Since $S$ was arbitrary,
  $\norm{\A\x}_2^2\geq\lambda_{\min}(\A_S^\top
  \A_S^{})\norm{\x}_2^2\geq\varepsilon\norm{\x}_2^2$ for all $\x$ with
  support $S\subseteq\{1,\dots,n\}$, $\abs{S}\leq k$. Moreover, the
  encoding length of $\alpha$, and therefore that of $\varepsilon$, is
  clearly polynomially bounded by the encoding length of $\A$, which
  completes the proof.
\end{IEEEproof}

\begin{theorem}\label{thm:RIPcertifyNP}
  Given a matrix $\A\in\Q^{m\times n}$, a positive integer $k$, and
  some constant $\delta_k\in (0,1)$, it is \NP-hard to decide
  whether $\A$ satisfies the RIP of order $k$ with
  constant~$\delta_k$.
\end{theorem}

\begin{IEEEproof}
  Consider a matrix $\tilde{\A}$ as in Lemma~\ref{lem:scaleA}, so we
  can again focus on the lower inequality of the RIP. Clearly, if
  $\tilde{\A}$ has a circuit of size at most $k$, $\tilde{\A}$ cannot
  satisfy the RIP of order $k$ with any given $\delta_k\in(0,1)$,
  since in this case,
  $\norm{\tilde{\A}\tilde{\x}}_2^2=0<(1-\delta_k)\norm{\tilde{\x}}_2^2$
  for some $\tilde{\x}$ with $1\leq\norm{\tilde{\x}}_0\leq k$. Moreover,
  Lemma~\ref{lem:RIPlowerBound} shows that if $\tilde{\A}$ has no
  circuit of size at most $k$, $\tilde{\A}$ satisfies the RIP of order
  $k$ with constant $1-\tilde{\varepsilon}\in(0,1)$, where
  $\tilde{\varepsilon}$ has size polynomially bounded by $k$ and that
  of~$\tilde{\A}$, cf.~\eqref{eq:RIPeps}.  By Theorem~\ref{thm:ripNP},
  it is \coNP-complete to decide whether there exists a constant
  $\delta_k<1$ such that $\tilde{\A}$ satisfies the RIP of a given
  order $k$ with this constant. But as seen above, such a constant
  exists if and only if $\tilde{\A}$ satisfies the RIP of $k$ with the
  constant $1-\tilde{\varepsilon}$, too. Thus, deciding whether the
  RIP holds for a given matrix, order, \emph{and} constant, is
  (\textsf{co})\NP-hard.
\end{IEEEproof}

\begin{remark}\label{rem:RIPcertifyCoNPcontainmentOpen}
  It is an open question whether the problem in
  Theorem~\ref{thm:RIPcertifyNP} is in \coNP.
\end{remark}

\begin{remark}
  Clearly, Theorem~\ref{thm:RIPcertifyNP} leads to another easy proof
  for Corollary~\ref{cor:ripNP01NPhard} (and
  Corollary~\ref{cor:asymmRIPlowerNP} below): computing the (lower
  asymmetric) RIC would also decide the RIP certification problem. On
  the other hand, our proof of Theorem~\ref{thm:RIPcertifyNP}
  essentially reduces the RIP certification problem to the setting of
  Theorem~\ref{thm:ripNP}, which therefore can be seen as the core RIP
  hardness result (by establishing the direct link to spark
  computations); see also Remark~\ref{rem:whyThm2} below.
\end{remark}

\subsection{Asymmetric restricted isometry constants}\label{subsect:asymmRIP}

It has been remarked in~\cite{BCT11} that the symmetric nature of the
RIP can be overly restrictive. In particular, the influence of the
upper inequality in~\eqref{eq:deltak} is often stronger, although the
lower inequality is more important in the context of sparse
recovery. For instance, the often stated condition
$\underline{\delta}_{2k} < 1$ for uniqueness of~$k$-sparse solutions
(see, e.g., \cite{C08}) should in fact read
$\underline{\delta}_{2k}^{L} < 1$, where
\begin{equation}\label{eq:deltakL}
  \underline{\delta}_{k}^{L}\define \min_{\delta\geq 0}\,\delta\quad\st\quad (1-\delta)\norm{\x}_2^2\leq\norm{\A\x}_2^2\quad\forall \x: \norm{\x}_0\leq k
\end{equation}
is the lower \emph{asymmetric} restricted isometry
constant~\cite{BCT11} (see also \cite{FL09}). Correspondingly, the
upper asymmetric RIC is 
\begin{equation}\label{eq:deltakU}
  \underline{\delta}_{k}^{U}\define \min_{\delta\geq 0}\,\delta\quad\st\quad (1+\delta)\norm{\x}_2^2\geq\norm{\A\x}_2^2\quad\forall \x: \norm{\x}_0\leq k.
\end{equation}
The central argument in the proof of Theorem~\ref{thm:ripNP} in fact
shows the following:
\begin{corollary}\label{cor:asymmRIPlowerNP}
  Given a matrix $\A\in\Q^{m\times n}$ and a positive integer $k$, it
  is \NP-hard to compute~$\underline{\delta}_k^L$.
\end{corollary}

Moreover, the next result settles the computational complexity of
computing the upper asymmetric RIC $\underline{\delta}^U_k$.
\begin{theorem}\label{thm:asymmRIPupperNP}
  Given a matrix $\A\in\Q^{m\times n}$, a positive integer $k\leq m$
  and a parameter $\delta>0$, it is \NP-hard in the strong sense
  to decide whether $\underline{\delta}^U_k<\delta$, even in the
  square case $m = n$. Consequently, it is strongly \NP-hard to
  compute $\underline{\delta}_k^U$.
\end{theorem}
To prove this theorem, we need some auxiliary results.
\begin{lemma}\label{lem:graphSpectra}
  Let $G=(V,E)$ be a simple undirected graph with $n=\abs{V}$ and let
  $\A_G$ be its $n\times n$ adjacency matrix, i.e., $(\A_G)_{ij}=1$ if
  $\{i,j\}\in E$ and $0$ otherwise. Denote by $K_n$ the complete graph with $n$ vertices.
  \begin{enumerate}
  \item If $G=K_n$, i.e., $G$ is a clique, then $\A_G$ has eigenvalues
    $-1$ and $n-1$ with respective multiplicities $n-1$ and $1$.
  \item Removing an edge from $G$ does not increase
    $\lambda_{\max}(\A_G)$. In fact, if $G$ is connected, this
    strictly decreases $\lambda_{\max}(\A_G)$.
  \item If $G=K_n\setminus e$, i.e., a clique with one edge
    removed, then the largest eigenvalue of $\A_G$ is
    $(n-3+\sqrt{n^2+2n-7})/2$.
  \end{enumerate}
\end{lemma}
\begin{IEEEproof}
  The first two statements can be found in, or deduced easily from,
  \cite[Ch.~1.4.1 and Prop.~3.1.1]{BH12}, respectively. The third
  result is a special case of \cite[Theorem~1]{F88}.
\end{IEEEproof}
\begin{remark}\label{rem:eigSepa}
  Lemma~\ref{lem:graphSpectra} shows that, in a graph $G=(V,E)$ with
  $\abs{V}=n$, the largest eigenvalue of the adjacency matrix of any
  induced subgraph with $k\in\{2,\dots,n\}$ vertices is either $k-1$
  (if and only if the subgraph is a $k$-clique) or at most
  $(k-3+\sqrt{k^2+2k-7})/2$.
\end{remark}
\begin{proposition}\label{prop:sparsePCAcoNP}
  Given a matrix $\H\in\Q^{n\times n}$, a positive integer $k\leq n$
  and a parameter $\lambda>0$, it is \coNP-complete in the strong
  sense to decide whether $\lambda_{\max}^{(k)}<\lambda$, where
  \begin{equation}\label{eq:PCAprob}
    \lambda_{\max}^{(k)} ~\define~ \max\{\;\x^\top \H \x\suchthat \norm{\x}_2^2=1,~\norm{\x}_0\leq k\;\} 
    ~=~ \max\{\;\lambda_{\max}(\H_{SS})\suchthat S\subseteq\{1,\dots,n\},~\abs{S}\leq k\;\}.
  \end{equation}
  Consequently, 
  solving the \emph{sparse principal component analysis (Sparse PCA)} problem
  is strongly \NP-hard.
\end{proposition}
\begin{IEEEproof}
  We reduce from the $k$-Clique Problem. Let $G=(V,E)$ be a simple
  undirected graph with $n$ vertices (w.l.o.g., $n\geq 2$). From $G$,
  construct its $n\times n$ adjacency matrix $\H\define\A_G$. By the previous
  Lemma, $G$ contains a $k$-clique if and only if
  $\lambda_{\max}^{(k)}(\H)=k-1\eqqcolon\lambda$. (Note that
  $\lambda_{\max}^{(k)}\leq k-1$ always holds by construction.) Hence,
  the Sparse PCA decision problem has a negative answer for the
  instance $(\H,k,\lambda)$ if and only if the $k$-Clique Problem has
  a positive answer. Since the latter is \NP-complete in the strong
  sense, and because all numbers appearing in the constructed Sparse
  PCA instance, and their encoding lengths, are polynomially bounded
  by $n$, the former is strongly (\textsf{co})\NP-hard.

  Moreover, consider a ``no'' instance of the Sparse PCA decision
  problem. Then, as we just saw, there is a $k$-clique $S$ in $G$, and
  it is easily verified that $\hat{\x}$ with $\hat{\x}_i=1/\sqrt{k}$
  for $i\in S$, and zeros everywhere else, achieves $\hat{\x}^\top \H
  \hat{x}=\lambda_{\max}^{(k)}(\H)=\lambda$. Scaling~\eqref{eq:PCAprob}
  by~$k$, we see that \emph{equivalently},
  \[
  k\cdot\lambda_{\max}^{(k)} ~=~ \max\{\;\x^\top \H\x\suchthat
  \norm{\x}_2^2=k,~\norm{\x}_0\leq k\;\} ~=~ k\lambda ~=~ k^2-k.
  \]
  Thus, a \emph{rational} certificate for the ``no'' answer is given
  by the vector $\tilde{\x}$ with $\tilde{\x}_i=1$ for $i\in S$, and
  zeros everywhere else. 
  Since we can clearly check all constraints on $\tilde{\x}$ (from the
  scaled problem) and that $\tilde{\x}^\top \H\tilde{\x}=k\lambda$ in
  polynomial time, the Sparse PCA decision problem is contained in
  \coNP.
\end{IEEEproof}
\begin{remark}\label{rem:sparsePCA}
  The Sparse PCA problem (see, e.g., \cite{LT13,dAEG11,dAEGJL07}) is
  often mentioned to be (\NP-)hard, but we could not locate a rigorous
  proof of this fact. In~\cite[Section~6]{BR12}, the authors sketch a
  reduction from the $k$-Clique Problem but do not give the details;
  the central spectral argument mentioned there, however, is exactly
  what we exploit in the above proof.
\end{remark}
We will extend the proof of Proposition~\ref{prop:sparsePCAcoNP} to show
Theorem~\ref{thm:asymmRIPupperNP} by suitably approximating the
Cholesky decomposition of a matrix very similar to the adjacency matrix; 
the following technical result will be useful for this extension.
\begin{lemma}\label{lem:CholeskyBounds}
  Let $\A_G$ be the adjacency matrix of a simple undirected graph
  $G=(V,E)$ with $n$ vertices, and let $\H\define \A_G+n^2 \I$, where
  $\I$ is the $n\times n$ identity matrix. Then $\H$ has a unique
  Cholesky factorization $\H=\L\D\L^\top$ with diagonal matrix
  $\D=\diag(d_1,\dots,d_n)\in\Q^{n\times n}$ and unit lower triangular
  matrix $\L=(\ell_{ij})\in\Q^{n\times n}$, whose respective entries
  fulfill $d_i\in[(n^4-2n+2)/n^2,n^2]$ and
  $\ell_{ij}\in[(2-n^2-2n)/(n^4-2n+2),(2n-2)/(n^4-2n+2)]$ for all
  $i,j\in\{1,\dots,n\}$, $i\neq j$.
\end{lemma}
\begin{IEEEproof}
  With deg$(v)$ denoting the degree of a vertex $v$ of $G$,
  and $\lambda_i(\A_G)$, $i=1,\dots,n$, the eigenvalues of $\A_G$,
  \[
  \norm{\A_G}_2 ~=~ \max_{1\leq i\leq n}\,\abs{\lambda_i(\A_G)} ~\leq~ \sqrt{n}\norm{\A_G}_{\infty} ~<~ n\cdot\max_{v\in V}\,\text{deg}(v) ~<~ n^2.
  \]
  Thus, it is easy to see that $\H=\A_G+n^2\I$ is (symmetric) positive
  definite; in particular, the eigenvalues of $\H$ obey
  $\lambda_i(\H)=\lambda_i(\A_G)+n^2$ for all~$i$. Then, $\H$ has a
  unique Cholesky factorization $\H=\L\D\L^\top$, and
  $\L\in\Q^{n\times n}$ (unit lower triangular) and $\D\in\Q^{n\times
    n}$ (diagonal) can be obtained by Gaussian elimination; see, e.g.,
  \cite[Section~4.9.2]{GMW91} or \cite[Section~4.2.3]{GvL96}.

  Let $\H^{(0)}\define\H$ and let $\H^{(k)}=(h^{(k)}_{ij})$ be the
  matrix obtained from $\H$ after $k$ iterations of (symmetric)
  Gaussian elimination. There are $n-1$ such iterations, and each
  matrix $\H^{(k)}$ has block structure with $\diag(d_1,\dots,d_k)$ in
  the upper left part and a symmetric matrix in the lower right block.
  We show by induction that for all $k=1,\dots,n$, and all
  $i,j\in\{0,\dots,n-k\}$,
  \begin{equation}\label{eq:CholInd}
    h^{(k-1)}_{k+i,k+i} \in \left[ n^2-\frac{2(k-1)}{n^2},n^2\right]\qquad\text{and (for $i\neq j$)}\qquad h^{(k-1)}_{k+i,k+j}\in\left[-\frac{2(k-1)}{n^2},1+\frac{2(k-1)}{n^2}\right].
  \end{equation}
  Clearly, by construction of $\H$, $h^{(0)}_{ii}=n^2$ and
  $h^{(0)}_{ij}\in\{0,1\}$ for all $i,j\in\{1,\dots,n\}$, $i\neq j$,
  so \eqref{eq:CholInd} holds true for $k=1$. Suppose
  \eqref{eq:CholInd} holds for some $k\in\{1,\dots,n-1\}$, i.e.,
  throughout the first $k-1$ iterations of the Gaussian elimination
  process. Performing the $k$-th iteration, we obtain
  $h^{(k)}_{ik}=h^{(k)}_{ki}=0$ for all $i>k$,
  $h^{(k)}_{kk}=h^{(k-1)}_{kk}$, and
  \begin{equation}\label{eq:chol_gen}
    h^{(k)}_{k+i,k+j} ~=~
    h^{(k-1)}_{k+i,k+j}-\frac{h^{(k-1)}_{k+i,k}}{h^{(k-1)}_{kk}}\cdot
    h^{(k-1)}_{k,k+j}\qquad\text{for all }i,j\in\{1,\dots,n-k\}.
  \end{equation}
  Thus, in particular, by symmetry of $\H^{(k-1)}$, for any $i\in\{1,\dots,n-k\}$,
  \begin{equation}\label{eq:chol_diag}
    h^{(k)}_{k+i,k+i} ~=~ h^{(k-1)}_{k+i,k+i}-\frac{h^{(k-1)}_{k+i,k}}{h^{(k-1)}_{kk}}\cdot h^{(k-1)}_{k,k+i} 
    ~=~ h^{(k-1)}_{k+i,k+i}-\frac{\big(h^{(k-1)}_{k+i,k}\big)^2}{h^{(k-1)}_{kk}}.
  \end{equation}
  Applying the induction hypothesis~\eqref{eq:CholInd}
  to~\eqref{eq:chol_diag} yields the first interval inclusion (note
  that $h^{(k-1)}_{kk}>0$, so $h^{(k)}_{k+i,k+i}\leq
  h^{(k-1)}_{k+i,k+i}$):
  \begin{equation}\label{eq:chol_diag_bounds}
  n^2 ~\geq~ h^{(k)}_{k+i,k+i} ~\geq~ n^2-\frac{2(k-1)}{n^2}-\frac{\big(1+\tfrac{2(k-1)}{n^2}\big)^2}{\big(n^2-\tfrac{2(k-1)}{n^2}\big)} ~=~ n^2-\frac{2(k-1)}{n^2}-\underbrace{\frac{n^2+4(k-1)+\tfrac{4(k-1)^2}{n^2}}{n^4-2k+2}}_{\leq~ 2/n^2} ~\geq~ n^2-\frac{2k}{n^2}.
  \end{equation}
  Similarly, for the off-diagonal entries $h^{(k)}_{k+i,k+j}$ with
  $i,j\in\{1,\dots,n-k\}$, $i\neq j$, from \eqref{eq:chol_gen}
  and~\eqref{eq:CholInd} we obtain
  \[
  h^{(k)}_{k+i,k+j} ~\leq~
  1+\frac{2(k-1)}{n^2}-\frac{\big(-\tfrac{2(k-1)}{n^2}\big)\big(1+\tfrac{2(k-1)}{n^2}\big)}{\big(n^2-\tfrac{2(k-1)}{n^2}\big)}
  ~=~ 1+\frac{2(k-1)}{n^2}+\underbrace{\frac{2(k-1)+\tfrac{4(k-1)^2}{n^2}}{n^4-2k+2}}_{\leq~ 2/n^2}
  ~\leq~ 1+\frac{2k}{n^2}
  \]
  and (compare with~\eqref{eq:chol_diag_bounds})
  \[
  h^{(k)}_{k+i,k+j} ~\geq~ -\frac{2(k-1)}{n^2}
  -\frac{\big(1+\tfrac{2(k-1)}{n^2}\big)^2}{\big(n^2-\tfrac{2(k-1)}{n^2}\big)}
  ~\geq~ \frac{2-2k}{n^2}-\frac{2}{n^2}
  ~=~ -\frac{2k}{n^2},
  \]
  which shows the second interval inclusion and concludes the induction.

  The statement of the Lemma now follows from observing that
  $d_k=h^{(k-1)}_{kk}$ for all $k=1,\dots,n$, and because the
  entries in the lower triangular part of $\L$, i.e., $\ell_{ij}$ with
  $i>j$, $j=1,\dots,n-1$, contain precisely the negated elimination
  coefficients (from the $j$-th iteration, respectively), whence
  \[
  \frac{2-n^2-2n}{n^4-2n+2}~=~-\frac{\big(1+\tfrac{2(n-1)}{n^2}\big)}{\big(n^2-\tfrac{2(n-1)}{n^2}\big)} ~\leq~ \ell_{ij} ~=~ -\frac{h^{(j-1)}_{ij}}{h^{(j-1)}_{jj}} 
    ~=~-\frac{h^{(j-1)}_{ij}}{d_j} ~\leq~ -\frac{\big(-\tfrac{2(n-1)}{n^2}\big)}{\big(n^2-\tfrac{2(n-1)}{n^2}\big)} ~=~ \frac{2n-2}{n^4-2n+2}.
  \]
  (By construction, $\ell_{ii}=1$ and $\ell_{ij}=0$ for all
  $i=1,\dots,n$, $j>i$.) This concludes the proof.
\end{IEEEproof}
\begin{IEEEproof}[Proof of Theorem~\ref{thm:asymmRIPupperNP}]
  Given an instance $(G,k)$ for the $k$-Clique Problem, we construct
  $\H=\A_G+n^2\I$ from the graph's adjacency matrix $\A_G$ (w.l.o.g.,
  $n\geq 2$). From the proof of Proposition~\ref{prop:sparsePCAcoNP},
  recall that $G$ has no $k$-clique if and only if
  $\lambda_{\max}^{(k)}(\A_G)<k-1$, or equivalently
  $\lambda_{\max}^{(k)}(\H)<n^2+k-1$ (cf. the beginning of the proof
  of Lemma~\ref{lem:CholeskyBounds}). Let $\D$ and $\L$ be the
  Cholesky factors of $\H$, i.e., $\H=\L\D\L^\top$. Letting
  $\D^{1/2}\define\diag(\sqrt{d_1},\dots,\sqrt{d_n})$, observe that
  the upper asymmetric RIC for the matrix $\A'\define\D^{1/2}\L^\top$
  can be written as
  \[
  \underline{\delta}^U_k(\A') ~=~ \max\{\;\x^\top \L\D^{1/2}\D^{1/2}\L^\top\x\suchthat\norm{\x}_2^2=1,~\norm{\x}_0\leq k\;\}-1 ~=~ \lambda_{\max}^{(k)}(\H)-1.
  \]
  Consequently, $G$ has a $k$-clique $S$ if and only if $\H$ has a
  $k\times k$ submatrix $\H_{SS}$ with largest eigenvalue $n^2+k-1$,
  i.e., $\underline{\delta}^U_k(\A')=n^2+k-2$. Moreover, by
  Lemma~\ref{lem:graphSpectra}, $\lambda_{\max}(\H_{TT})\leq
  n^2+(k-3+\sqrt{k^2+2k-7})/2<n^2 +k-1$ 
  for any incomplete induced subgraph of $G$ with vertex set $T$,
  $\abs{T}=k$. However, while $\L$ and $\D$ are rational, $\D^{1/2}$
  can contain \emph{irrational} entries, so we cannot directly use
  $\A'$ as the input matrix for the upper asymmetric RIC decision
  problem.  The remainder of this proof shows that we can replace
  $\D^{1/2}$ by a rational approximation to within an accuracy that
  still allows us to distinguish between eigenvalues associated to
  $k$-cliques and those belonging to incomplete induced subgraphs.
    
  Let us consider the rational approximation obtained by truncating
  after the $p$-th decimal number (we will specify $p$ later), i.e.,
  let $\tilde{\D}^{1/2}\define\diag(r_1,\dots,r_n)$ with $r_i\define
  \lfloor 10^p\cdot\sqrt{d_i}\rfloor/10^p$. Thus, $\sqrt{d_i}-r_i\leq
  10^{-p}$ for all $i$ by construction, and in particular, since
  $d_1=n^2$ (see Lemma~\ref{lem:CholeskyBounds}),
  $r_1=n$. Consequently,
  \[
  \norm{\D^{1/2}-\tilde{\D}^{1/2}}_2 ~=~ \max_{1\leq i\leq n}\{\abs{\sqrt{d_i}-r_i}\} ~=~ \max_{2\leq i\leq n}\{ \sqrt{d_i}-r_i\} ~\leq~ 10^{-p}.
  \]
  Denoting $\tilde{\D}\define \tilde{\D}^{1/2}\tilde{\D}^{1/2}$ and using $d_i\leq n^2$ (by Lemma~\ref{lem:CholeskyBounds}), we obtain
  \[
  \norm{\D-\tilde{\D}}_2 ~=~ \max_{2\leq i\leq n}\{ d_i-r_i^2 \} ~=~ \max_{2\leq i\leq n}\{(\sqrt{d_i}+r_i)(\sqrt{d_i}-r_i)\} ~\leq~ 10^{-p}\cdot 2\cdot\max\{\sqrt{d_i},r_i\} ~\leq~ 2\cdot 10^{-p}\cdot n.
  \]
  Let $\tilde{\H}\define \L\tilde{\D}\L^\top$ and note that, for all
  $i,j=1,\dots,n$, we have
  $h_{ij}=\sum_{q=1}^{n}d_q\ell_{iq}\ell_{jq}$ and
  $\tilde{h}_{ij}=\sum_{q=1}^n r_q^2 \ell_{iq}\ell_{jq}$. Since
  $\abs{\ell_{ij}}\leq 1$ for all $i,j$ (by
  Lemma~\ref{lem:CholeskyBounds}), and by symmetry of $\H$ and
  $\tilde{\H}$, it follows that
  \[
  \abs{(\H-\tilde{\H})_{ij}} ~=~ \abs{(\H-\tilde{\H})_{ji}} ~=~ \Big\lvert\sum_{q=1}^{n}(d_q-r_q^2)\ell_{iq}\ell_{jq}\Big\rvert ~\leq~ \sum_{q=1}^{n}\abs{d_q-r_q^2}\,\abs{\ell_{iq}}\,\abs{\ell_{jq}} ~\leq~ \sum_{q=1}^n (d_q-r_q^2)~\leq~ 2\cdot 10^{-p}\cdot n^2.
  \]
  Thus, we have $\tilde{\H}=\H+\tilde{\E}$, where
  $\abs{\tilde{\E}_{ij}}\leq 2\cdot 10^{-p}\cdot n^2$ and $\tilde{\E}$
  is also symmetric. Note that, for any $S\subseteq\{1,\dots,n\}$,
  \[
  \lambda_{\max}(\tilde{\E}_{SS}) ~\leq~ \lambda_{\max}(\tilde{\E}) ~=~ \norm{\tilde{\E}}_2 ~\leq~ \sqrt{\norm{\tilde{\E}}_1\cdot\norm{\tilde{\E}}_\infty} ~\leq~
  \sqrt{(n\cdot 2\cdot 10^{-p}\cdot n^2)(n\cdot 2\cdot 10^{-p}\cdot n^2)} ~=~ 2\cdot 10^{-p}\cdot n^3.
  \]
  Therefore (cf., e.g.,~\cite[Corollary~8.1.6]{GvL96}), we have for all $S\subseteq\{1,\dots,n\}$ that        
  \[
  \abs{\lambda_i(\H_{SS})-\lambda_i(\tilde{\H}_{SS})} ~\leq~ \norm{\tilde{\E}_{SS}}_2 ~=~ \lambda_{\max}(\tilde{\E}_{SS}) ~\leq~ 2\cdot 10^{-p}\cdot n^3\qquad\text{for all }i=1,\dots,\abs{S}.
  \]
  Consequently, if $S$ is a $k$-clique, we have
  \begin{equation}\label{eq:lambdaClique}
  \lambda_{\max}(\tilde{\H}_{SS}) ~\geq~ n^2+k-1-2\cdot 10^{-p}\cdot n^3,
  \end{equation}
  whereas for any $T\subseteq\{1,\dots,n\}$ with $\abs{T}=k$ which
  induces no clique, it holds that
  \begin{equation}\label{eq:lambdaNoClique}
  \lambda_{\max}(\tilde{\H}_{TT}) ~\leq~ n^2+(k-3+\sqrt{k^2+2k-7})/2 + 2\cdot 10^{-p}\cdot n^3. 
  \end{equation}
  
  Now fix $p\define 1+\lceil 4\log_{10}(n)\rceil$, and let an instance
  for the upper asymmetric RIC decision problem be given by
  $\A\define\tilde{\D}^{1/2}\L^\top$ (where $\tilde{\D}^{1/2}$ is
  computed from the Cholesky factors $\L$ and $\D$ of $\H=\A_G+n^2\I$
  using this precision parameter~$p$), $\delta\define n^2+k-2-2\cdot
  10^{-p}\cdot n^3$, and $k$.
  
  If $G$ has a $k$-clique, then by~\eqref{eq:lambdaClique},
  $\underline{\delta}^U_k(\A)\geq\delta$, and if not, 
  $\underline{\delta}^U_k(\A)\leq n^2+(k-3+\sqrt{k^2+2k-7})/2 + 2\cdot 10^{-p}\cdot n^3-1$ 
  by~\eqref{eq:lambdaNoClique}. In fact, our choice of $p$ implies
  \[
  p ~>~ \log_{10}(8)+4\log_{10}(n) \quad\Rightarrow\quad 10^p ~>~ 8n^4 ~>~ \frac{8n^3}{k+1-\sqrt{k^2+2k-7}},
  \]
  from which we can derive that
  \[
  n^2+k-1-2\cdot 10^{-p} \cdot n^3 ~>~ n^2 + \frac{k-3+\sqrt{k^2+2k-7}}{2}+2\cdot 10^{-p}\cdot n^3,
  \]
  which shows that $\underline{\delta}^U_k(\A)<\delta$ if no
  $k$-clique is contained in $G$. Therefore, $G$ has a $k$-clique if
  and only if $\underline{\delta}^U_k(\A)\geq\delta$, i.e., the upper
  asymmetric RIC decision problem under consideration has a negative
  answer.

  Clearly, all computations in the above reduction can be performed in
  polynomial time. To see that the encoding length~$\langle\A\rangle$
  of $\A$ is in fact polynomially bounded by $n$, note that
  $h_{ij}\in\{0,1,n^2\}$ for all $i,j$, whence
  $\langle\H\rangle\in\O(\poly(n))$. Since Gaussian elimination can be
  implemented to lead only to a polynomial growth of the encoding
  lengths~\cite{GroLS93}, it follows that
  $\langle\L\rangle$,$\langle\D\rangle\in\O(\poly(\langle\H\rangle))=\O(\poly(n))$.
  In particular, the entries of $\tilde{\D}^{1/2}$ then also have
  encoding length polynomially bounded by $n$, by construction of the
  rational approximation. This shows that indeed $\langle
  a_{ij}\rangle\in\O(\poly(n))$ for all $i,j$. Furthermore, all the
  numerical values $a_{ij}$ are also polynomially bounded by $n$ (in
  fact, $\abs{a_{ij}}\leq n$, by Lemma~\ref{lem:CholeskyBounds} and
  the construction of $\tilde{\D}^{1/2}$). Moreover, $0<\delta<n^2+n$
  and, clearly, its encoding length $\langle\delta\rangle$ is bounded
  polynomially by $n$ as well.

  Thus, since the Clique Problem is well-known to be \NP-complete in
  the strong sense, and because our polynomial reduction in fact
  preserves boundedness of the numbers within $\O(\poly(n))$, the
  upper asymmetric RIC decision problem is (\textsf{co})\NP-hard in
  the strong sense.  
  This completes the proof of
  Theorem~\ref{thm:asymmRIPupperNP}.
\end{IEEEproof}


In fact, observe that for the matrix $\A$ constructed in the proof of
Theorem~\ref{thm:asymmRIPupperNP}, the upper asymmetric RIC is always
larger than the lower asymmetric RIC, whence the former therefore
coincides with the (symmetric) RIC $\underline{\delta}_k$ of
$\A$. Thus, the following result holds true, which slightly
strengthens Corollary~\ref{cor:ripNP01NPhard}.
\begin{corollary}\label{cor:RICstronglyNP}
  Given a matrix $A\in\Q^{m\times n}$ and a positive integer $k$, it
  is \NP-hard in the strong sense to compute the RIC $\underline{\delta}_k$, even
  in the square case $m=n$.
\end{corollary}

\begin{remark}\label{rem:whyThm2}
  Theorem~\ref{thm:ripNP} (which yielded
  Corollary~\ref{cor:ripNP01NPhard}) is of interest in its own right,
  as it reveals, for instance, the intrinsic relation between the
  $k$-sparse solution uniqueness conditions $2k<\spark(\A)$ and
  $\underline{\delta}_{2k}<1$ (or $\underline{\delta}_{2k}^{L}<1$,
  respectively, see Corollary~\ref{cor:asymmRIPlowerNP});
  consequently, verifying uniqueness via these conditions is \NP-hard
  because deciding whether $\spark(\A)\leq k$ is.
\end{remark}

\begin{remark}
  We can easily extend the construction from the proof of
  Theorem~\ref{thm:asymmRIPupperNP} to cover the non-square case $m<n$
  as well. The idea is as follows: For instance, let $\B$ be the
  matrix from Example~\ref{ex:spark}, and replace $\A$ in the above
  proof by the block diagonal matrix $\hat{\A}$ which has $\A$ in the
  first block and $\B$ in the second. This matrix has dimensions
  $(n+3)\times (n+4)$, and since the eigenvalues of $\B^\top \B$ are
  contained in $[0,5)$, it is easy to see that the relation between
  eigenvalues of symmetrically chosen submatrices and cliques is the
  same for $\hat{\A}^\top\hat{\A}$ as for $\A^\top \A$ itself (note
  that $n^2+(k-3+\sqrt{k^2+2k-7})/2-2\cdot 10^{-p}\cdot n^3\geq 5$
  whenever $n\geq 3$, $k\geq 2$, which can of course be assumed
  without loss of generality).
  A similar idea is exploited in the proof of~\cite[Theorem~6]{KZ12}.
\end{remark}

\begin{remark}
  The rational certificate for the Sparse PCA problem (see the proof
  of Proposition~\ref{prop:sparsePCAcoNP}) cannot be employed to
  verify the ``no'' answer of the upper asymmetric RIC decision
  problem from Theorem~\ref{thm:asymmRIPupperNP} in polynomial time;
  containment in \coNP\ of this latter problem thus remains an open
  question.
\end{remark}

\section{\NP-hardness of computing the nullspace constant}\label{sect:nspNP}

Recall that~$\A \in \Q^{m \times n}$ satisfies the NSP of order~$k$
with constant~$\alpha_k$ if~\eqref{eq:NSPineq} holds for all $\x \in
\R^n$ with $\A\x = 0$. As for the RIP, one is typically interested in
the smallest constant~$\underline{\alpha}_k$, the \emph{nullspace
  constant (NSC)}, such that the NSP of order~$k$ holds with this
constant. Thus, the NSC is given~by
\begin{equation}\label{eq:alphak}
  \underline{\alpha}_k \define \min\, \alpha \quad\st\quad \norm{\x}_{k,1} \leq \alpha
  \norm{\x}_1 \quad\text{for all }\x\text{ with }\A\x = 0,
\end{equation}
or equivalently,
\begin{equation}\label{eq:alphak2}
  \underline{\alpha}_k \define \max ~\norm{\x_S}_1 \quad\st\quad \A\x = 0,\; \norm{\x}_1=1,\; S \subseteq \{1, \dots, n\},\; \abs{S} \leq k.
\end{equation}
Sparse recovery by \eqref{prob:P1} is ensured if and only if
$\underline{\alpha}_k<1/2$. However, the following results show that
computing $\underline{\alpha}_k$ is a challenging problem.

\begin{theorem}\label{thm:nspNP}
  Given a matrix $\A \in \Q^{m \times n}$ and a positive integer $k$,
  the problem to decide whether~$\A$ satisfies the NSP of order~$k$
  with some constant $\alpha_k<1$ is \coNP-complete.
\end{theorem}

\begin{IEEEproof}
  First of all, note that the NSP~\eqref{eq:NSPineq} is equivalent to
  the condition that
  \begin{equation}\label{eq:NSPeqiv}
    \norm{\x_S}_1 ~\leq~ \alpha_k\, \norm{\x}_1
  \end{equation}
  holds for all $S \subseteq \{1, \dots, n\}$, $\abs{S} \leq k$, and
  $\x \in \R^n$ with $\A\x = 0$. Clearly,~\eqref{eq:NSPineq}
  and~\eqref{eq:NSPeqiv} always hold for some $\alpha_k \in [0,1]$.

  We first show that the problem is in \coNP. To certify the ``no''
  answer, it suffices to consider a vector~$\tilde{\x}$ with $\A
  \tilde{\x} = 0$ and a set~$\emptyset\neq S \subseteq \{1, \dots,
  n\}$ with $\abs{S} \leq k$ that tightly satisfy~\eqref{eq:NSPeqiv}
  for~$\alpha_k = 1$. This implies that~$S$ contains the support
  of~$\tilde{\x}$. Thus, $1\leq\norm{\tilde{\x}}_0 \leq k$. Clearly,
  since $\tilde{\x}$ is contained in the nullspace of~$\A$, it can be
  assumed to be rational with encoding length polynomially bounded by
  that of~$\A$. This shows that the ``no'' answer can be certified
  in polynomial time.

  To show hardness, we claim that the matrix~$\A$ has a circuit of
  size at most~$k$ if and only if there does \emph{not} exist
  any~$\alpha_k < 1$ such that~\eqref{eq:NSPineq} holds. Since the
  former problem is \NP-complete by Theorem~\ref{thm:circuitNP}, this
  completes the proof.
  
  Assume $\A$ has a circuit of size at most $k$. Then there exists a
  vector~$\x$ in the nullspace of~$\A$ with $1\leq\norm{\x}_0 \leq
  k$. It follows that $\norm{\x}_{k,1} = \norm{\x}_1$. Thus,
  \eqref{eq:NSPineq} implies $\alpha_k \geq 1$. Since, trivially,
  $\alpha_k \leq 1$, we conclude that $\alpha_k = 1$.
  
  Conversely, assume that there exists no $\alpha_k < 1$ such
  that~\eqref{eq:NSPineq} holds for $\A$ and~$k$. This implies that
  there is a vector~$\x$ with $\A\x = 0$ and $1\leq\norm{\x}_0 \leq k$
  such that $\norm{\x}_{k,1} = \norm{\x}_1$, because otherwise
  $\alpha_k < 1$ would be possible. But this means that the support
  of~$\x$ contains a circuit of $\A$ of size at most~$k$, which shows
  the claim.
\end{IEEEproof}

We immediately obtain the following.

\begin{corollary}\label{cor:nspNPhard}
  Given a matrix $\A \in \Q^{m \times n}$ and a positive integer~$k$,
  it is \NP-hard to compute the nullspace
  constant~$\underline{\alpha}_k$.
\end{corollary}

\section{Concluding remarks}\label{sect:conclusion}

The results of this paper show that it is \coNP-complete to answer the
following questions in the case $\gamma = 1$: ``Given a matrix $\A$
and a positive integer $k$, does the RIP or NSP hold with some
constant $< \gamma$?'' It is important to note that our results do not
imply the hardness for every \emph{fixed} constant~$\gamma < 1$. (Note
that such questions become solvable in time $\O(n^{\text{poly}(k)})$
if~$k$ is fixed.) For instance, Theorem~\ref{thm:RIPcertifyNP} asserts
that it is (\textsf{co})\NP-hard to certify the RIP for given $\A$,
$k$ and $\delta_k \in (0,1)$ \emph{in general}. The actual $\delta_k$
appearing in the proof (as in the one independently derived
in~\cite{BDMS12}), however, is very close to $1$ and thus far from
values of $\gamma$ that could yield recovery guarantees. Similarly,
the NSP guarantees \emph{$\ell_0$-$\ell_1$-equivalence} for
$\underline{\alpha}_k < 1/2$, while we proved \NP-hardness only for
deciding whether $\underline{\alpha}_k < 1$ (which implies that
\emph{computing} $\underline{\alpha}_k$ is \NP-hard). The complexity
of these related questions remains open. Nevertheless, our results
provide a justification for investigating general approximation
algorithms (which compute bounds on $\underline{\delta}_k$ or
$\underline{\alpha}_k$), as done in \cite{dAEG11,dAEGJL07,JN11,LB08},
instead of searching for exact polynomial time algorithms.

Instead of the intractable RIP, NSP, or spark, the weaker but
efficiently computable \emph{mutual coherence} \cite{DH01} is
sometimes used. It can be shown that the mutual coherence yields
bounds on the RIC, NSC, and the spark; see, for instance,
\cite{GN03,BDE09}. Thus, imposing certain conditions involving the
mutual coherence of a matrix can yield uniqueness and recoverability
(by basis pursuit or other heuristics), see,
e.g.,~\cite{BDE09}. However, the sparsity levels for which the mutual
coherence can guarantee recoverability are quite often too small to be
of practical use. This emphasizes the importance of other concepts.

An interesting question for future research is whether it is hard to
approximate the constants associated with the RIP or NSP in polynomial
time to within some factor, or whether spark and NSC computations are
also strongly \NP-hard (as is RIC computation, see
Corollary~\ref{cor:RICstronglyNP}). A first step in this direction was
taken in~\cite{KZ12}, where inapproximability of RIP parameters is
shown under certain less common complexity assumptions (interestingly,
also making use of Cholesky decompositions of certain matrices related
to a type of adjacency matrix for random graphs), see
also~\cite{KZ11,BR13}; strong inapproximability results
for~\eqref{prob:P0} appear in~\cite{AK98}.

Moreover, (\textsf{co})\NP-hardness does not necessarily exclude the
possibility of \emph{practically} efficient algorithms. So far, to the
best of our knowledge, the focus has been laid largely on relaxations
or heuristics. In~\cite{KZ11,KZ12}, it was shown that one may
sometimes do better than exhaustive search to certify the RIP, making
use of the nondecreasing monotonicity of $\underline{\delta}_k$ with
growing~$k$. A ``sandwiching'' procedure to compute the NSC
$\underline{\alpha}_k$ (exactly) was very recently proposed
in~\cite{CX13} and empirically demonstrated to be faster than brute
force. However, neither method can \emph{guarantee} a running time
improvement with respect to simple enumeration. Moreover,
Corollary~\ref{cor:RICstronglyNP} shows that, in general, no
pseudo-polynomial time algorithm (i.e., a method with running time
polynomially bounded by the input size and the largest occuring
numerical value) can exist for computing the RIC
$\underline{\delta}_k$, unless \P$=$\NP. More work on exact algorithms
could shed more light on the behavior of the RIP and NSP.

\section*{Acknowledgment} 
We thank the anonymous referees for their constructive comments and
for spotting an error in an earlier version of this paper.

\bibliographystyle{IEEEtran}
\bibliography{IEEEabrv}

%
%
%
%
%
\begin{IEEEbiographynophoto}{Andreas M. Tillmann}
  graduated in financial and business mathematics
  (Dipl.-Math. Oec. degree) from the TU Braunschweig, Germany, in
  2009. As a doctoral candidate at the TU Braunschweig, he worked as a
  Teaching Assistant from 2009 to 2011, and was a Researcher on the
  project ``SPEAR -- Sparse Exact and Approximate Recovery'', funded
  by the Deutsche Forschungsgemeinschaft (DFG), from 03/2011 to
  06/2013. In 07/2012 he has joined the Research Group Optimization at
  TU Darmstadt, Germany; in 07/2013 he resumed a TA position
  there. His research interests are mainly in compressed sensing,
  focussing on exact $\ell_0$-minimization, basis pursuit, and issues
  related to recovery conditions.
\end{IEEEbiographynophoto}

\begin{IEEEbiographynophoto}{Marc E. Pfetsch}
  obtained a Diploma degree in mathematics from the University of
  Heidelberg, Germany, in 1997. He received a PhD degree in
  mathematics in 2002 and the Habilitation degree in 2008 from TU
  Berlin, Germany. From 2008 to 2012 he was a Full Professor for
  Mathematical Optimization at TU Braunschweig, Germany. Since 04/2012
  he has been a Full Professor for Discrete Optimization at TU
  Darmstadt, Germany. His research interests are mostly in discrete
  optimization, in particular symmetry in integer programs, compressed
  sensing, and algorithms for mixed integer programs.
\end{IEEEbiographynophoto}



\vfill


\end{document}